\documentclass[pdflatex,sn-mathphys-num]{sn-jnl}

\usepackage{graphicx}%
\usepackage{multirow}%
\usepackage{amsmath,amssymb,amsfonts}%
\usepackage{amsthm}%
\usepackage{mathrsfs}%
\usepackage[title]{appendix}%
\usepackage{xcolor}%
\usepackage{textcomp}%
\usepackage{manyfoot}%
\usepackage{booktabs}%
\usepackage{algorithm}%
\usepackage{algorithmicx}%
\usepackage{algpseudocode}%
\usepackage{listings}%
\usepackage{comment}

\theoremstyle{thmstyleone}%
\newtheorem{theorem}{Theorem}
\newtheorem{proposition}[theorem]{Proposition}%

\theoremstyle{thmstyletwo}%
\newtheorem{remark}{Remark}%

\theoremstyle{thmstylethree}%
\newtheorem{definition}{Definition}%

\usepackage{caption}
\def\R{\mathbb R}

\def\int{{\rm int\,}}

\def\eps{\varepsilon}

\theoremstyle{plain}
\theoremstyle{theorem}
\numberwithin{theorem}{section}
\newtheorem{corollary}[theorem]{Corollary}
\newtheorem{lemma}[theorem]{Lemma}
\theoremstyle{definition}

\begin{document}

\title[Article Title]{A fixed-time stable forward-backward  dynamical  system for solving generalized monotone inclusions}


\author[3]{\fnm{Nam} \sur{V. Tran}}\email{namtv@hcmute.edu.vn}

\author[1,2,3]{\fnm{Le} \sur{T.T. Hai}}\email{hailtt@hcmute.edu.vn}

\author[3]{\fnm{Truong} \sur{V. An}}\email{antv@hcmute.edu.vn}

\author*[4]{\fnm{Phan} \sur{ T.  Vuong}}\email{t.v.phan@soton.ac.uk}

{\color{red}
\affil[1]{\orgdiv{Faculty of Mathematics and Computer Science}, \orgname{University of Science}, \orgaddress{\street{Nguyen Van Cu}, \city{District 5}, \postcode{70000}, \state{Ho Chi Minh City}, \country{Vietnam}}}

\affil[2]{\orgname{Vietnam National University}, \orgaddress{\street{Linh Trung}, \city{Thu Duc}, \postcode{70000}, \state{Ho Chi Minh City}, \country{Vietnam}}}
}

\affil[3]{\orgdiv{Faculty of Applied Sciences}, \orgname{{\color{black} HCMC University of Technology and Education}, \orgaddress{\street{Vo Van Ngan}, \city{Thu Duc}, \postcode{70000}, \state{Ho Chi Minh City}, \country{Vietnam}}}}

\affil*[4]{\orgdiv{School of Mathematical Sciences}, \orgname{University of Southampton}, \orgaddress{\street{University Street}, \city{Southampton}, \postcode{SO17 1BJ}, \state{Hampshire}, \country{UK}}}


\abstract{ We propose a  forward-backward splitting dynamical system for solving inclusion problems of the form $0\in A(x)+B(x)$ in Hilbert spaces, where $A$ is a maximal operator and $B$ is a single-valued operator. Involved operators are assumed to satisfy  a generalized monotonicity condition, which is weaker than the standard monotone assumptions.  Under  mild conditions on parameters, we establish the fixed-time stability of the proposed dynamical system. In addition, we consider an explicit forward Euler discretization of the  dynamical system leading to a new forward backward algorithm for which we present the convergence analysis. Applications to other optimization problems such as Constrained Optimization Problems (COPs), Mixed Variational Inequalities (MVIs), and Variational Inequalities (VIs) are presented 
and some numerical examples are given to illustrate the theoretical results.

}




\maketitle

	\section{Introduction} 
	Let $H$ be a real  finite-dimensional Hilbert space endowed with an inner product and its induced norm denoted by
	$\langle \cdot,\cdot \rangle$ and $\|\cdot\|$, respectively. 
	In this article we are interested in the following  inclusion problem: Find \(x^*\in H\) such that 
	\begin{equation}\label{inc1}
		0\in A(x^*) + B(x^*)
	\end{equation}
	where $A:H \longrightarrow 2^H$ is a multi-valued operator and $B:H \longrightarrow H$ is a single-valued operator. This problem has an important role in many fields such as optimization problems, variational inequalities, equilibrium problems,  saddle point problems, 
	 Nash equilibrium problem in noncooperative games,  fixed point problems, 
	and others; see, for instance, \cite{BlumOettli94, Cav} and references quoted therein. 
	
For example, let $f: H \longrightarrow (-\infty, +\infty], f\not\equiv +\infty $ be a subdifferential and convex function,  and $C$ be a closed, convex subset of $H$, then the constrained optimization problem (COP): $\min_{x\in C} f(x)$ is equivalent to the monotone inclusion 
$0\in  \partial f (x)+N_C(x)$, where \(N_C(x)\) is the \emph{normal cone} to \(C\) at \(x\) (see, for instance, \cite{ROCK}). 
	
	Also, a  variational inequality problem: find \(x^*\in C\) such that \(\langle F(x^*), x-x^*\rangle~\geq~0, \quad \forall x\in C\), where \(C\) is a closed, convex subset of \(H\), can be expressed as an inclusion problem $0\in A(x)$ where 
	\begin{equation*}
		A(x)=\begin{cases} 
			F(x)+N_C(x) &\mbox{ if } x\in C\\
			\emptyset &\mbox{ if } x\neq  C.
		\end{cases}
	\end{equation*} 
	
	
	One of the most common methods for solving problem \eqref{inc1} is the forward-backward splitting algorithm, introduced by Passty \cite{PAS} and Lions and Mercier \cite{MER}. Afterward, many modified versions of this algorithm are developed by numerous authors \cite{CSE, MAL, TSENG}.
	
	In recent years, dynamical systems have gained popularity due to their high performance and low computational cost in addressing various problems such as inclusion problems (IPs), mixed variational inequalities (MVIs), variational inequalities (VIs), and constrained optimization problems (COPs) \cite{AAS, BC_SICON, 22,   BCJMAA, BotCV,    18,  21,   27, 26, 20, JU4, LIU, 29,  VS20, ZHU}. For example, Bo\c t et al. \cite{22} investigated the global exponential stability of the forward-backward-forward dynamical system proposed in \cite{24}, while in \cite{BC_SICON}, the authors applied second-order forward-backward dynamical systems to inclusion problems. Additionally, \cite{BOT3} presented the strong convergence of continuous Newton-like inertial dynamics with Tikhonov regularization for monotone inclusions. Despite the widespread application of dynamical systems, much of the theoretical research on these methods {\color{blue} have} focused on their asymptotic stability \cite{27, 26, 29} or exponential stability \cite{22, 18,  21, 20, JU, JU4,  LIU,   V, VS20, ZHU}.

	Incorporating the concept of finite-time convergence, as proposed in \cite{BHAT}, several dynamical systems \cite{JU6, LI} have been recently introduced, exhibiting finite-time convergence, where the settling time relies on initial states and may escalate indefinitely with the increasing deviation of the initial state from an equilibrium point. However, in numerous practical scenarios like robotics and vehicle monitoring networks, obtaining accurate initial states {\color{blue} are} often challenging, if not impossible, beforehand. Consequently, determining the settling time in advance becomes impractical. To address this limitation, the concept of fixed-time convergence was introduced, wherein the convergence time can be upper-bounded by a constant independent of initial states. This concept was initially proposed by Polyakov in his seminal work \cite{polyakov} and has since been further developed by numerous researchers, see, for example, \cite{ Garg21, JU3, SAN, WANG, WANG2, WANG3, ZUO}.
	
	The principles of fixed-time convergence and finite-time convergence have also found application in various optimization and control contexts (see, for instance, multi-agent control problems \cite{ZUO} and sparse optimization problems \cite{GARG4,  HE, YU}). For instance, Chen and Ren \cite{CHEN} introduced a sign projected gradient method exhibiting finite-time convergence to address a category of convex optimization problems. Garg and Panagou \cite{GAR} proposed a gradient flow with fixed-time convergence for solving unconstrained optimization problems. Romero and Benosman \cite{ROMEO} explored four dynamical systems demonstrating finite-time convergence, bridging the gap between the q-rescaled gradient flow \cite{WIB} and the normalized gradient flow \cite{COR}. Furthermore, Garg et al. \cite{Garg21} devised a fixed-time convergent proximal dynamical system for addressing MVIs. Subsequently, Garg and Baranwal \cite{GARG4}  applied it to tackle sparse optimization problems with a non-smooth regularizer. More recently, Ju et al. \cite{JU3}  extended the proximal dynamical system \cite{JU4} and the forward-backward-forward projection dynamical system \cite{22} to a fixed-time converging forward–backward–forward proximal dynamical system for solving MVIs. X. Ju et al. applied the proximal dynamical system to study fixed-time convergence  \cite{JU5} and finite-time convergence  \cite{JU6} for equilibrium problems (EPs). 
	
	The literature review presented above highlights the consideration of dynamical system with finite-time or fixed-time convergence for solving various problems such as MVIs and VIs \cite{Garg21, JU3}, EPs \cite{JU6, JU5}, COPs \cite{JU6}, linear programming problems \cite{SAN}, unconstrained optimization problems (UOPs) \cite{CHEN, COR, GAR, ROMEO, WIB}, and sparse signal recovery problems \cite{ GARG4, HE, YU}. However, their application in solving inclusion problems has not been investigated. Specifically, the development of fixed-time converging forward-backward splitting dynamical systems for inclusion problems remains unexplored, which serves as the motivation for this study.

	In this paper, we continue this research direction by considering a forward-backward splitting first-order dynamical system associated with a fixed point reformulation (see e.g. \cite{AM,Bot,BSV,V0}), for which we obtain the fixed-time stability. 
	In addition, we consider a discrete version of the proposed dynamical system, which leads to a forward-backward method with relaxation. We establish the linear convergence of the iterative sequence generated by this algorithm to the unique solution of the inclusion problem.
	
	The monotonicity assumption in inclusion problems is standard and challenging to abandon, as even common results may prove invalid without this hypothesis. For instance, if $(A+B)$ lacks strong monotonicity, the inclusion $0\in (A+B)(x)$ may have no solution. In this paper, we will adopt a generalized notion of monotonicity, allowing the modulus of monotone operators to be negative. This broadened perspective enables us to address a wider range of operators beyond those traditionally considered monotone operators.

	The contributions of this article, in contrast to existing relevant research, can be summarized as follows.
	\begin{enumerate}
		\item This work introduces a novel forward-backward splitting dynamical system with fixed-time convergence for addressing inclusion problems.  To the best of our knowledge, this is the first study to propose a fixed-time convergent forward-backward splitting dynamical system for solving inclusion problems. 
  In contrast to existing forward-backward splitting dynamical system methods for solving inclusions, where the settling time for convergence to the solution cannot be explicitly determined, this work provides an explicit upper bound for the settling time required to reach the solution.
		
		\item Some applications to COPs, MVIs, VIs are given. By treating these problems as special cases of inclusion problems, we explore the fixed-time stability of those problems and provide upper bounds for settling time. 
		
		\item Unlike the monotonicity or strong monotonicity assumptions commonly imposed on operators in prior works, this article introduces a requirement of generalized monotonicity known as $\mu$-monotonicity, wherein the modulus of the operator may take negative values. Monotonicity has emerged as a standard assumption in the study of inclusion problems. Generalized monotonicity enables us to broaden the scope of operators considered in these inclusion problems

	\end{enumerate}

	The remaining sections of the article are structured as follows. Section \ref{Preliminaries} begins by revisiting some fundamental definitions and concepts,  forward-backward splitting dynamical systems, finite-time, and fixed-time stabilities, along with the introduction of several technical lemmas. In Section \ref{sec3}, we present the main results, including the newly proposed fixed time forward-backward splitting dynamical system model and its analysis. 
	Section \ref{sec4} introduces a condition that is adequate for attaining a consistent discretization, meaning a discretization that converges with a fixed number of time steps, for dynamical systems represented by a broad category of differential inclusions. As a specific example, it demonstrates that the forward-Euler discretization of the modified proximal dynamical system is indeed consistent.
	Section \ref{sec5} explores the applications of the forward-backward dynamical system in solving COPs, MVIs, {\color{blue} and} VIs. 

 To demonstrate the effectiveness and advantages of the proposed finite-time converging proximal dynamical system,  Section \ref{num} offers some numerical examples.  
	Finally, in Section \ref{sec6}, we provide some concluding remarks.
	
\section{Preliminaries} \label{Preliminaries}
	
	In this section, we recall some well-known definitions useful in the sequel.

	\subsection{Some notions on convex analysis}

	Let $g: H  \to \mathbb {R} \cup \{+\infty\}$ be a proper, convex and lower semicontinuous (l.s.c.) function.  We call $g$ \emph{subdifferentiable} at $x$ if the set
	$$
	\partial g(x)= \{ u \in H: \, g(y) \ge  g(x)+ \left\langle u,y - x \right\rangle  \, \forall y  \in H \}
	$$
	is nonempty. Then, $\partial g(x)$ is called the \emph{subdifferential} of $g$ at $x$, and vector $u \in \partial g(x)$ is called a \emph{subgradient} of $g$ at $x$. The function $g$ is subdifferentiable on $H$ if it is subdifferentiable at each point of $H$. 
		Note that if the function $g$ is convex, l.s.c. and has the full domain, it is continuous on the whole space \cite[Corollary 8.30]{BauschkeCombettes}. In this case, $g$ is subdifferentiable on $H$ \cite[Proposition 16.14]{BauschkeCombettes}. 
		In addition, let $f$ and $g$ be proper, convex, l.s.c functions such that $\text{dom} f \cap \text{int dom} g \not= \emptyset$ or $\text{dom}(g) = H$, 
		here $\text{dom} f = \{ x \in H, f(x) < +\infty\}$ denotes  the domain of $f$, 
		then 
		$\partial (f+g) = \partial f + \partial g$ \cite[Corollary 16.38]{BauschkeCombettes}. 
		
		
		Let $C$ be a closed, convex subset of $H$. The normal cone $N_C$ to $C$ at a
		point $x\in C$ is defined by 
		$$ N_C(x)=\left\{w\in H:\left\langle w,x-y\right\rangle \ge 0, \forall y\in C\right\}, $$
		and $N_C(x) = \emptyset$ if $x \not\in C$. 
		The indicator function of $C$ is defined as 
		$i_C(x) =0$ if $x\in C$ and $i_C(x) =+\infty$ otherwise. Then, we have 
		$\partial i_C(x) = N_C(x)$ for all $x\in H$.\\
		
		For every $x\in H$, the
		metric projection $P_C(x)$ of $x$ onto $C$ is defined by
		$$
		P_C (x)=\arg\min\left\{\left\|y-x\right\|:y\in C\right\}.
		$$
		Since $C$ is nonempty, closed, and convex, $P_C (x)$ exists and is unique. 
  
Finally, let $f: H \to \overline{\mathbb{R}}$. The proximity operator of parameter $\lambda >0$ of a function $f$ at $x \in H$ is defined by 
	\begin{eqnarray*}
		\mbox{prox}_{\lambda f}(x)=\mbox{argmin}_{y\in H}\left\{f(y)+\dfrac{1}{2\lambda}\|x-y\|^2\right\}\quad \forall x\in H.
	\end{eqnarray*}
For more details
		as well as for unexplained terminologies and notations we refer to \cite{BauschkeCombettes}.

	\subsection{Monotone operators}

	Let  $A: H \longrightarrow 2^H$ be a given operator. 	Recall that the domain of \(A\) is denoted by \(\mbox{dom }A\) and defined by 
	$$\mbox{dom }A=\{x\in H: A(x)\neq \emptyset\}.$$ We define  $G(A)=\{(x, u)\in H\times H: u\in A(x)\}$ the \emph{graph} of \(A\). 
	
	\begin{definition}
		Let  $A: H \longrightarrow 2^H$ be a given operator. The operator \(A\) is said to be  \emph{ $\mu_A$ - monotone} if   $\langle x-y, u-v\rangle \geq \mu_A\|x-y\|^2$ for all $x, y\in H, u\in A(x), v\in A(y)$ and some $\mu_A \in \R$. 	
	\end{definition}

	\begin{definition}
		A $\mu_A$-monotone operator \(A\) is said to be \emph{maximal} if its graph is not properly contained in the graph of any $\mu_A$-monotone operator \(A': H\longrightarrow 2^H\) \cite{Minh}. 
	\end{definition}
	
	\begin{remark} Observe that $\mu_A$ in the above definition can be negative. If $\mu_A=0$, $\mu_A$-monotonicity reduces to the classical monotonicity. If $\mu_A>0$, A (maximal) $\mu_A$-monotonicity  degenerates to (maximal) strong monotonicity. Finally, if $\mu_A<0$, $A$ is called \emph{weakly-monotone}. 
	\end{remark}
	
	One example of  a maximal monotone  operator is the subdifferential $\partial f$ of a lower semicontinuous convex function $f: H \longrightarrow (-\infty, +\infty], f\not\equiv +\infty $ (see, e.g. \cite{ROCK}). 
	
	\begin{definition}
		An operator $B:H\longrightarrow H$ is said to be 
		\begin{itemize}
			\item \emph{Lipschitz continuous} with constant \(L\geq 0\) if $\|B(x)-B(y)\|\leq L\|x-y\|$ for all \(x,y\in H\).
			\item  \emph{$\beta$-cocoercive} {\color{blue} with constant $\beta\geq 0$} if $\langle B(x)-B(y), x-y\rangle \geq \beta\|B(x)-B(y)\|^2$ for all $x, y\in H$. 
		\end{itemize}
	\end{definition}
	
By the Cauchy inequality, it is clear that if an operator is $\beta$-cocoercive, then it is $\frac{1}{\beta}$ Lipschitz continuous.\\	
	We recall here the resolvent of an operator $A: H \longrightarrow 2^H$, which is a useful tool in studying inclusion problems. The resolvent  with the parameter $\lambda$  of the operator $A$ is defined as follows 
	$$J_{\lambda A} =(Id+\lambda A)^{-1},$$
	where $Id$ is denoted by the \emph{identity mapping}. 

Without monotonicity, the resolvent of a maximal operator may not be single-valued.  The lemma below shows that the singleton property of the resolvent is maintained for generalized monotone operators with the appropriate values of the parameter. It also shows that the resolvent is cocoercive. This result will be applied several times in the sequel. 
 
	\begin{lemma}\label{GENMON} \cite{Minh} Let $A : H \longrightarrow 2^H$ be a $\mu_A$-monotone operator and let $\lambda  > 0$ be such that $1+\lambda \mu_A>0$. Then, the following statements hold
		\begin{enumerate}
			\item $J_{\lambda A}$ is single-valued;
			\item $J_{\lambda A}$ is $(1 + \lambda \mu_A)$-cocoercive;
			\item ${\color{black}\mbox{\rm dom}} J_{\lambda A} =H$  if and only if $A$ is maximal $\mu_A$-monotone
		\end{enumerate}
	\end{lemma}
	
	For more details on monotone operators, their application to optimization problems, and properties of their resolvent we refer the readers to \cite{AVR, BauschkeCombettes, Minh}.
	
	We denote by \(\text{\rm Fix}(A)\) the \emph{set of fixed points} and by $\mbox{zer} A$ the \emph{set of zero points} of an operator \(A\).
	The following lemma allows us to characterize zero points of the operator $A+B$ as fixed points of an appropriate associated operator.
	
	\begin{lemma}   \label{lem fixed point}
		Let $A: H\longrightarrow 2^H$  be  $\mu_A$-monotone and $B: H\longrightarrow H$. Let \(\lambda >0\) be such that $1+\lambda \mu_A>0$. Then $x^*\in H$ is a zero point of inclusion problem \eqref{inc1} if and only if $J_{\lambda A}(x^*-\lambda B(x^*))=x^*$, i.e., if and only if $x^*\in {\color{black}\text{\rm Fix}}(J_{\lambda A}\circ (Id-\lambda B))$.
	\end{lemma}
	\begin{proof}
Since \(A\) is \(\mu_A\)-monotone and $\lambda >0, 1+\lambda \mu_A>0$, by Lemma \ref{GENMON}, $J_{\lambda A}$ is a singleton. Then $0\in A(x)+B(x)$ if and only if $-B(x)\in A(x)$ which is equivalent to $x-\lambda B(x)\in x+\lambda A(x)$ or equivalent to $x\in \big((Id+\lambda A)^{-1}\circ (Id-\lambda B)\big)(x)= J_{\lambda A}\circ (Id-\lambda B)(x)$. The proof is completed. 
	\end{proof}

	
	\subsection{Equilibrium points and stability}
	To establish the main results of the paper, we need to recall the stability concepts of an equilibrium point of the general dynamical system
	\begin{equation}\label{JJS1}
		\dot{x}(t) = T(x(t)), \quad t \ge 0,
	\end{equation}
	where $T$ is a continuous mapping from $H$ to $H$ and $x:[0, +\infty)\to H$.
	\begin{definition} \cite{Pappaladro02}
		\begin{description}
			\item[{\rm (a)}] A point $x^*\in H$ is an equilibrium point for (\ref{JJS1}) if $T(x^*)=0$;
			\item[{\rm (b)}] An equilibrium point $x^*$ of (\ref{JJS1}) is stable if, for any $\epsilon>0$, there exists  $r>0$ such that, for every $x_0 \in B(x^*, r) $, the solution $x(t)$ of the dynamical system with $x(0)=x_0$ exists and is contained in  $B(x^*, \epsilon)$ for all $t>0$, where $B(x^*, r )$ denotes the open ball with center $x^*$ and radius $r$;
			\item[{\rm (c)}]  A stable equilibrium point $x^*$ of (\ref{JJS1}) is asymptotically stable if there exists
			$r>0$ such that, for every solution $x(t)$ of (\ref{JJS1}) with $x(0) \in B(x^*, r )$, one has
			$$ 
			\lim_{t \to +\infty} x(t)=x^*;
			$$
			\item[{\rm (d)}]  An equilibrium point $x^*$ of (\ref{JJS1}) is exponentially stable if there exist
			$r>0$  and constants $\kappa>0$ and $\theta>0$ such that, for every solution
			$x(t)$ of (\ref{JJS1}) with $x(0) \in B(x^*, r )$, one has
			\begin{equation} \label{exponentialstability}
				\| x(t) -x^* \| \leq \kappa\, \| x(0)-x^*\|\, e^{- \theta t} \quad \forall t \geq 0.
			\end{equation}
			Furthermore, $x^*$ is globally exponentially stable if (\ref{exponentialstability}) holds true for all solutions $x(t)$ of (\ref{JJS1}).
		\end{description}
	\end{definition}
	\begin{definition}
		The equilibrium point \(x^*\) is said to be
		\begin{enumerate}[(a)]
			\item {\emph{Finite-time stable}} if it is stable in the sense of Lyapunov, and there exists a neighborhood \(B(r, x^*)\) of \(x^*\) and a settling-time function \(T: B(r, x^*)\setminus \{x^*\}\rightarrow (0, \infty)\) such that for any \(x(0)\in  B(r, x^*)\setminus\{x^*\}\), the solution of \eqref{JJS1} satisfies \(x(t)\in B(r, x^*)\setminus \{x^*\}\) for all \(t\in [0, T(x(0)))\) and \(\lim_{t\to T(x(0))}x(t)=x^*\).
			\item \emph{Globally finite-time stable} {\color{blue} if } it is finite-time stable with \(B(r, x^*)=\R^n\).
			\item \emph{ Fixed-time stable} if it is globally finite-time stable, and the settling-time function satisfies
			\[\sup_{x(0)\in \R^n}T(x(0))<\infty.\]
		\end{enumerate}
	\end{definition}
	
 Polyakov obtained the following results for fixed-stable time \cite{polyakov}:
	\begin{theorem} \label{lm1}	(Lyapunov condition for fixed-time stability). Suppose that there exists 	a continuously differentiable function $V : D\rightarrow \R$, where $D \subseteq H$  is a neighborhood of the	equilibrium point \(x^*\) for \eqref{JJS1} such that
		\begin{equation*}
			V(x^*) = 0, V(x) > 0
		\end{equation*}
		for all $x \in  D \setminus \{x^*\}$ and
		\begin{equation} \label{inq of Lyapunov} 
			\dot{V}(x) \leq  -(p_1V(x)^{\alpha_1} +p_2V(x)^{\alpha_2})^{\alpha_3}
		\end{equation}
		for all $x \in D \setminus \{x^*\}$ with $p_1,p_2,\alpha_1,\alpha_2,\alpha_3 > 0$ such that $\alpha_1\alpha_3 <1$ and $\alpha_2\alpha_3 >1.$  Then, the equilibrium point \(x^*\) of \eqref{JJS1} is fixed-time stable such that
		\begin{equation}\label{time1}
			T(x(0)) \leq \dfrac{1}{p_1^{\alpha_3}(1-\alpha_1\alpha_3)} +\dfrac{1}{p_2^{\alpha_3}(\alpha_2\alpha_3-1)}
		\end{equation}
		for any $x(0) \in H$.  
		Moreover, take $\alpha_1=\left(1-\frac{1}{2\xi}\right), \alpha_2=\left(1+\frac{1}{2\xi}\right),$ with $\xi>1$ and $\alpha_3=1$ in \eqref{inq of Lyapunov}. Then the equilibrium point of \eqref{JJS1} is fixed-time stable with the settling time 
		\begin{equation}\label{time2}
			T\leq T_{\max} =\dfrac{\pi\xi}{\sqrt{p_1p_2}}.
		\end{equation}

		In addition, if the function $V$ is radially unbounded (i.e., $||x||\to \infty \Rightarrow V(x)\to \infty$) and $D =H$, then the equilibrium point \(x^*\) of \eqref{JJS1} is globally fixed-time stable.
	\end{theorem}
	
	It was shown in \cite{PARS} that the settling time \eqref{time1} gives a conservative convergence times estimate, meanwhile settling time \eqref{time2} presents a less conservative one for dynamical system \eqref{JJS1}. 
	We will apply this fundamental result to analyze the fixed-time stability of solutions of inclusion  problem \eqref{inc1} via appropriate dynamical systems in the next section.

	\section{ Forward-backward dynamical systems and fixed-time stability} 	\label{sec3}
	In this section, we propose a novel forward-backward splitting dynamical system model.  We then characterize the solutions of inclusion problems via equilibrium points of appropriate dynamical systems. Finally, the fixed-time convergence of the  dynamical system is established.
	
	\subsection{A nominal dynamical system for inclusion problems }
	
	For solving inclusion \eqref{inc1} we consider the nominal dynamical system 
	
	\begin{equation}\label{firstdynamicalsystem} 
		\dot{x}(t)=-\sigma \left(x(t)-J_{\lambda A}(x(t)-\lambda B(x(t)))\right), 
	\end{equation}
	where \(\sigma\) is a positive constant. 
 
 The first result in this section gives a characterization for solutions of inclusion {\color{black} problem} \eqref{inc1} via equilibrium points of the dynamical system \eqref{firstdynamicalsystem}. 
	\begin{proposition}\label{rem2} 
	 Let $A: H\longrightarrow 2^H$ be a $\mu_A$-monotone, let $B: H \longrightarrow H$ and  $\lambda >0$ be such that $1+\lambda \mu_A>0$. Then a point \(x^*\in H\) is a solution of inclusion \eqref{inc1} if and only if it is an equilibrium point of \eqref{firstdynamicalsystem}.
	\end{proposition}
		
\begin{proof}
It follows from  Lemma \ref{lem fixed point} and the definition of equilibrium points of the dynamical system \eqref{firstdynamicalsystem}. 
\end{proof}

In the sequel we will use the below assumption.

	\begin{itemize}
\item [{\bf(A)}]
		$A$ is maximal $\mu_A$-monotone, $B$ is $\mu_B$-monotone and $L$-Lipschitz continuous and the parameter $\lambda>0$ such that   
    $1+\lambda \mu_A>0 $ and $ 2(\mu_A + \mu_B) + \lambda \mu_A^2 - \lambda L^2 > 0. 
  $	
	\end{itemize}
	
The lemma below provides conditions for the existence of $\lambda>0$ satisfying assumption {\bf (A)}. 
\begin{lemma}\label{dktt lambda}
 Let $A: H\longrightarrow 2^H$ be maximal $\mu_A$- monotone and $B: H \longrightarrow H$ be  $\mu_B$-monotone and $L$-Lipschitz continuous. Then there exists a parameter $\lambda>0$ satisfying assumption {\bf (A)}  if and only if one of the following conditions on parameters $\mu_A, \mu_B, L$ is valid. 
  \begin{itemize}
      \item [{\bf (B1)}]  $\mu_A+\mu_B>0. $
        \item [{\bf (B2)}] $\mu_A+\mu_B=0$ and  $\mu_A^2>L^2$.
        \item [{\bf (B3)}] $\mu_A+\mu_B<0 $ and $\mu_A>L$.
  \end{itemize}
\end{lemma}
\begin{proof} Assume that there exists $\lambda>0$ satisfying assumption {\bf (A)}. There are only three cases that may occur:  $\mu_A+\mu_B>0, \mu_A+\mu_B=0$, and $\mu_A+\mu_B<0$.  If the first case occurs, {\bf(B1)} holds. 

When $\mu_A+\mu_B=0$, the last condition in assumption {\bf(A)} and the existence of $\lambda >0$ imply that   $\mu_A^2>L^2$. Hence {\bf (B2)} holds.

  As for the last case where $\mu_A+\mu_B<0 $, the last condition in assumption {\bf (A)} and the positiveness of $\lambda$ imply that $\mu_A^2>L^2$ which is equivalent to $\mu_A<-L$ or $\mu_A>L$.

If $\mu_A<-L$, there exists  $\lambda>0$ satisfying conditions in assumption {\bf (A)} if and only if there exists $\lambda$ such that $\dfrac{-2(\mu_A+\mu_B)}{\mu_A^2-L^2}<\lambda <-\dfrac{1}{\mu_A}$ hold. These inequalities hold for some $\lambda>0$ if and only if \begin{equation}\label{ineq for lambda}
    L^2+\mu_A^2+2\mu_A\mu_B<0
\end{equation} holds. This inequality implies that $0<\mu_B$.  Because of $\mu_B$-monotonicty and $L$-Lipschitz continuity of $B$, one has $\mu_B\leq L$. It follows that  $L^2+\mu_A^2+2\mu_A\mu_B\geq 0$, for all $\mu_A$. This means that the inequality \eqref{ineq for lambda} can not be satisfied and hence, there does not exist $\lambda>0$ satisfying assumption {\bf (A)} when $\mu_A+\mu_B<0, \mu_A<-L$.  Thus {\bf(B3)} holds.

We now prove for the sufficient condition of existence of $\lambda>0$ fulfilled {\bf(A)}. We assume that one of {\bf (Bi)}, $i=1,2,3$ holds. We will consider each case in turn. 

If $\mu_A+\mu_B>0$, the existence of $\lambda >0$ satisfying assumption {\bf (A)} is obviously. In particular, it holds for all $\lambda \in \left(0, \dfrac{2(\mu_A+\mu_B)}{L^2}\right)$.

If $\mu_A+\mu_B=0$ and  $\mu_A^2>L^2$, then the second condition on $\lambda$ in {\bf (A)} is satisfied for all $\lambda>0$. Also, the other condition on $\lambda$ is fulfilled for some  $\lambda>0$ and for any $\mu_A$. This assures the existence of $\lambda>0$. 

Finally, if $\mu_A+\mu_B<0 $ and $\mu_A>L$, the last condition in {\bf (A)} is fulfilled for $\lambda>\dfrac{-2(\mu_A+\mu_B)}{\mu_A^2-L^2}>0$ and the first condition on $\lambda$ is fulfilled for all $\lambda >0$. {\color{black} Thus} the existence of $\lambda$ satisfying the conditions in {\bf (A)} is assured. 
  \end{proof}

\begin{remark} From Lemma \eqref{dktt lambda}, note that condition {\bf (A)} does not imply the (strong) monotonicity of $A+B$. If  \(A\) is maximal monotone, (i.e. $\mu_A=0$) and \(B\) is strongly monotone (i.e. \(\mu_B>0\)),  then assumption {\bf (A)} holds for all $\lambda \in \left(0,\dfrac{2\mu_B}{L^2}\right)$.
\end{remark}

	The following proposition plays an important role in our convergence analysis of the proposed method.  
	
	\begin{proposition} \label{tr1}  Let \(A: H \longrightarrow 2^H\), $B:H \longrightarrow H$ and the parameter \(\lambda >0 \) be such that assumption {\bf (A)} holds.  	Then,  there exists $\delta\in (0,1)$  such that 
		\begin{enumerate}
			\item The operator $T:=J_{\lambda A}\circ (Id-\lambda B)$ is Lipchitz continuous with constant $\delta$.
   \item Fixed points of the operator $T$ exist and are uniquely determined. 
			\item $\|J_{\lambda A}(x-\lambda B(x))- x^*\|\leq \delta\|x-x^*\|$ for all $x^*\in {\color{black}\text{\rm Fix}}(T)$. 
		\end{enumerate}
	\end{proposition} 

	\begin{proof} (1) By Lemma \ref{GENMON}, and the choice of $\lambda$ we have that $J_{\lambda A}$ is single-valued and $(1+\lambda\mu_A)$-cocoercive. This implies that $J_{\lambda A}$ is $\dfrac{1}{1+\lambda \mu_A}$- Lipschitz continuous. We derive that  
		\begin{equation}\label{dg1}
			\begin{array}{lll}
				& \quad \|J_{\lambda A}(x-\lambda B(x))-J_{\lambda A}(y-\lambda B(y))\|^2  \\

				&	\leq \dfrac{1}{(1+\lambda \mu_A)^2}
    \|(x-\lambda B(x)) -(y-\lambda B(y))\|^2 
				\\
				&=\dfrac{1}{(1+\lambda \mu_A)^2}\left(\|x-y\|^2-2\lambda \langle x-y, B(x)-B(y)\rangle+\lambda^2 \|B(x)-B(y)\|^2\right)\\
				&\leq \dfrac{1}{(1+\lambda \mu_A)^2}\Big( \|x-y\|^2-2\mu_B \lambda \|x-y\|^2 +L^2\lambda^2\|x-y\|^2\Big)\\
				&=	\dfrac{1}{(1+\lambda \mu_A)^2}\Big(1-\lambda\left(2\mu_B -\lambda L^2\right)\Big)\|x-y\|^2 \quad \forall \ x, y\in H. 
			\end{array}		
		\end{equation} 
Note that under assumption {\bf (A}) 
it holds that 
		$$
		\dfrac{1-\lambda(2\mu_B -\lambda L^2)}{(1+\lambda \mu_A)^2} \in [0,1).
		$$   
		From the estimate \eqref{dg1} we conclude that $T$ is $\tau$-Lipschitz continuous with constant $$\tau=\sqrt{ \dfrac{1-\lambda(2\mu_B -\lambda L^2)}{(1+\lambda \mu_A)^2}}, \mbox{ where }\tau\in [0, 1).$$  Consequently, $T$ is $\delta$-Lipschitz continuous with some constant $\delta \in (0,1)$. 
  
(2) By Part (1),		$T$ is a contraction operator. Also, $A$ is maximal $\mu_A$-monotone, by Lemma \ref{GENMON}, $\mbox{dom } J_{\lambda A}=H$.  This implies that the fixed point $x^*$ of $T$  exists and is unique.
		
(3) It follows from Part (1) by noting that \(x^*\in {\color{black}\text{\rm Fix}}(T)\). 
	\end{proof}

The next result presents conditions for the existence and uniqueness of zero points of the sum of two generalized monotone operators.

 \begin{corollary}
Let \(A: H \longrightarrow 2^H\) be maximal $\mu_A$-monotone, and $B: H \longrightarrow H$  be $\mu_B$-monotone, $L$-Lipschitz continuous. Assume  that one of the following conditions holds
\begin{enumerate}
    \item $\mu_A+\mu_B>0, $
        \item  $\mu_A+\mu_B=0$ and  $\mu_A^2>L^2$,
        \item $\mu_A+\mu_B<0 $ and $\mu_A>L$.
\end{enumerate}
   
   Then solutions of inclusion problem \eqref{inc1} exist and are uniquely determined. 
 \end{corollary}
 \begin{proof}
   If one of the conditions is satisfied, by  Lemma \ref{dktt lambda}, there exists $\lambda>0$ satisfying assumption {\bf (A)}. Hence, Proposition \ref{tr1} can be invoked. Consequently, $J_{\lambda A}\circ (Id-\lambda B)$ has a unique fixed point.  It follows from Lemma \ref{lem fixed point} that the inclusion {\color{blue} problem} \eqref{inc1} has a unique solution. 
 \end{proof}

	\subsection{A novel dynamical system for fixed time stability of solutions of inclusion problems}
	In order to analyze the fixed-time stability of solutions to inclusion problems of the form \eqref{inc1} we introduce a modified dynamical system as follows
	\begin{equation}\label{newydynamicalsystem} 
		\dot{x}=-\varphi(x)(x-J_{\lambda A}(x-\lambda B(x))),
	\end{equation} 
	where 
	\begin{equation}\label{dnrho} 
		\varphi(x):=\begin{cases} 
			c_1\dfrac{1}{\|x-J_{\lambda A}(x-\lambda B(x))\|^{1-\kappa_1}}& +c_2\dfrac{1}{\|x-J_{\lambda A}(x-\lambda B(x))\|^{1-\kappa_2}} \\
			& \mbox{ if } x\in H\setminus {\rm Fix}(J_{\lambda A}\circ (Id-\lambda B))\\
			0 & \mbox{otherwise} 
		\end{cases}
	\end{equation} 
	with \(c_1, c_2>0, \kappa_1\in (0, 1)\) and \(\kappa_2>1\).
	
	The next result shows the relationship between the equilibrium points of the dynamical system above and those of \eqref{firstdynamicalsystem}.
	\begin{proposition}\label{pro3}
		A point $ x\in H$ is an equilibrium point of \eqref{newydynamicalsystem} if and only if it is an equilibrium point of \eqref{firstdynamicalsystem}. 
	\end{proposition}
	\begin{proof} It follows from the definition \eqref{dnrho} of $\varphi$ with noting that \(\varphi(x)=0\) if and only if \(x\in \text{\rm Fix}(J_{\lambda A}\circ (Id-\lambda B))\). 
	\end{proof}
Combining Proposition \ref{pro3} and Proposition \ref{rem2}  gives us a characterization for solutions of the inclusion {\color{blue} problem} \eqref{inc1} via equilibrium points of the dynamical system \eqref{newydynamicalsystem}. 

\begin{proposition}\label{pro4}
     Let \(A: H \longrightarrow 2^H\), $B:H \longrightarrow H$ and the parameter $\lambda>0$ be such that  assumption {\bf (A)} holds.  Then \(x^*\) is a solution of inclusion problem \eqref{inc1} if and only if it is an equilibrium point of the dynamical system \eqref{newydynamicalsystem}. 
\end{proposition}
\begin{proof}
    It follows from Proposition \ref{pro3} and Proposition \ref{rem2}. 
\end{proof}

\subsection{Fixed-time stability analysis} 
	We start this section by giving conditions so that solutions in the classical sense of the given dynamical system exist and are uniquely determined.
	
	\begin{lemma}\cite{Garg21}\label{prop.1}
		Let $T: H \rightarrow H$ be a locally Lipschitz continuous vector field such that 
		$$	T(\bar{x}) = 0 \mbox{ and } \langle x -\bar{x}, T(x) \rangle  > 0$$
		for all $x \in H\setminus\{\bar{x}\}$. Consider the following autonomous differential equation:
		\begin{equation}\label{system3}
			\dot{x}(t)=-\beta(x(t)) T(x(t)),
		\end{equation}
		where 
		\begin{equation*}
			\beta(x):=\begin{cases}
				c_1\dfrac{1}{\|T(x)\|^{1-\kappa_1}}+c_2\dfrac{1}{\|T(x))\|^{1-\kappa_2}} & \mbox{ if } T(x)\neq 0\\ 
				0 & \mbox{otherwise} 
			\end{cases}
		\end{equation*} 
		with \(c_1, c_2>0, \kappa_1\in (0, 1)\) and \(\kappa_2>1\). Then, with any given initial condition, the solution of \eqref{system3} exists in the classical sense and is uniquely determined for all $t\geq 0$. 
	\end{lemma}
	\begin{remark}\label{rm1} 
		Setting \(T(x):=x-J_{\lambda A}(x-\lambda B(x))\) for all \(x\in H\). Then this vector field satisfies the following property:
		\begin{equation*}
			\langle x-\overline{x}, T(x) \rangle > 0
		\end{equation*}
		for all $x\in H\setminus \{\overline{x}\}$, where \(\{\overline{x}\}= \rm Fix(J_{\lambda A}\circ (Id-\lambda B))\). 
		
		Indeed,  by Proposition \ref{tr1}, the set \({\color{black}\text{\rm Fix}}(J_{\lambda A}\circ (Id-\lambda B))\) is a singleton. Consequently,  the vector field in \eqref{firstdynamicalsystem} has a unique equilibrium point \(\overline{x}=x^*\).
		
		\noindent In addition, for all \(x\in H\setminus \{\overline{x}\},\) it holds that 
		\begin{equation}\label{eqn1} 
			\langle x-\overline{x}, x-J_{\lambda A}(x-\lambda B(x)) \rangle =\|x-\overline{x}\|^2+\langle x-\overline{x}, \overline{x}- J_{\lambda A} (x-\lambda B(x))\rangle.
		\end{equation}
		By the Cauchy-Schwarz inequality and Proposition \ref{tr1} with noting that \(\overline{x}=x^*\) one has 
		\begin{equation}\label{eqn2} 
			\langle x-\overline{x}, \overline{x}-J_{\lambda A} (x-\lambda B(x))\rangle \geq -\|x-\overline{x}\|\|\overline{x}-J_{\lambda A} (x-\lambda B(x))\| \geq -\delta \|x-\overline{x}\|^2 
		\end{equation} for some \(\delta\in (0,1)\). 
		Combining \eqref{eqn1} and \eqref{eqn2} one has 
		\begin{equation*}
			\langle x-\overline{x}, x-J_{\lambda A}(x-\lambda B(x)) \rangle \geq (1-\delta)\|x-\overline{x}\|^2>0 
		\end{equation*} for all \(x\in H\setminus \{\overline{x}\}\). 
	\end{remark}
	The following lemma is needed for the proof of the main theorem. 
	\begin{lemma}\cite{Garg21}\label{lm4} 
		For every \(\delta\in (0, 1)\), there is  \(\eps(\delta)=\dfrac{\log(\delta)}{\log\frac{1-\delta}{1+\delta}}>0\) satisfying that 
		\begin{equation}\label{eq11} 
			\left(\dfrac{1-\delta}{1+\delta}\right)^{1-\gamma}>\delta
		\end{equation}
		for any \(\gamma\in (1-\eps(\delta), 1)\). Moreover, \eqref{eq11} is valid for any \(\delta\in (0,1)\) and \(\gamma>1\). 
	\end{lemma}

	We are now ready for the main theorem of this section.
	
	\begin{theorem}\label{tr2}   Let $A: H\longrightarrow 2^H$, $B: H \longrightarrow H$ and \(\lambda \) be such that assumptions {\bf(A)} holds.  Then,  there exists \(\eps>0\) such that the solution \(x^*\in H\) of the inclusion problem \eqref{inc1} is a fixed-time stable equilibrium point of \eqref{firstdynamicalsystem} for any \(\kappa_1\in (1-\eps, 1)\cap(0, 1)\) and \(\kappa_2>1\) and the following time estimate holds:
		$$T(x(0))\leq T_{\max}=\dfrac{1}{p_1(1-\alpha_1)}+\dfrac{1}{p_2(\alpha_2-1)}$$
		for some $p_1>0, p_2>0, \alpha_1\in (0.5, 1), \alpha_2>1$. 
		
		In addition, if take $\alpha_1=1-\frac{1}{2\xi}, \alpha_2=1+\frac{1}{2\xi}$ with $\xi>1$   then the following time estimate holds:
		$$T(x(0))\leq T_{\max} =\dfrac{\pi \xi}{\sqrt{p_1p_2}}$$
		for some constants $p_1>0, p_2>0$ and $\xi>1$. 
	\end{theorem}
	\begin{proof} By Proposition \ref{tr1}, the operator \(J_{\lambda A}\circ(Id-\lambda B)\) is Lipschitz continuous with constant \(\delta\in (0, 1)\). It follows that the vector field at the right-hand side of \eqref{firstdynamicalsystem} is Lipschitz continuous on \(H\). By Remark \ref{rm1} the assumptions in Proposition \ref{prop.1} are satisfied. As a result, a solution of \eqref{newydynamicalsystem} exists and is uniquely determined for all forward times. 
		
		Consider the Lyapunov function \(V: H \rightarrow \R\) defined as follows:
		\[V(x):=\dfrac{1}{2}\|x-x^*\|^2.\]
		The time-derivative of the Lyapunov function \(V\) along the solution of \eqref{newydynamicalsystem}, starting from any \(x(0)\in H \setminus\{x^*\}\) with noting that $x^*\in \text{\rm Fix}(J_{\lambda A}(Id-\lambda B))$ being unique  yields:
		\begin{align}\label{eq12} 
			\dot{V}=&-\left\langle x-x^*, c_1\dfrac{x-J_{\lambda A}(x-\lambda B(x))}{\|x-J_{\lambda A}(x-\lambda B(x))\|^{1-\kappa_1}}+c_2\dfrac{x-J_{\lambda A}(x-\lambda B(x))}{\|x-J_{\lambda A}(x-\lambda B(x))\|^{1-\kappa_2}}\right\rangle \notag\\
			=&-\left\langle x-x^*, c_1\dfrac{x-x^*}{\|x-J_{\lambda A}(x-\lambda B(x))\|^{1-\kappa_1}}+c_2\dfrac{x-x^*}{\|x-J_{\lambda A}(x-\lambda B(x))\|^{1-\kappa_2}}\right\rangle \notag \\
			&-\left\langle x-x^*, c_1\dfrac{x^*-J_{\lambda A}(x-\lambda B(x))}{\|x-J_{\lambda A}(x-\lambda B(x))\|^{1-\kappa_1}}+c_2\dfrac{x^*-J_{\lambda A}(x-\lambda B(x))}{\|x-J_{\lambda A}(x-\lambda B(x))\|^{1-\kappa_2}}\right\rangle,
		\end{align}
  for all $x\in H\setminus \{x^*\}$. 
		By Cauchy-Schwarz inequality, we obtain an upper bound for the second term in the right hand side of \eqref{eq12} as follows
		\begin{align*}
			\dot{V}\leq & -\left(c_1\dfrac{\|x-x^*\|^2}{\|x-J_{\lambda A}(x-\lambda B(x))\|^{1-\kappa_1}}+c_2\dfrac{\|x-x^*\|^2}{\|x-J_{\lambda A}(x-\lambda B(x))\|^{1-\kappa_2}}\right) \notag \\
			&+ \left(c_1\dfrac{\|x-x^*\| \cdot \|x^*-J_{\lambda A}(x-\lambda B(x))\|}{\|x-J_{\lambda A}(x-\lambda B(x))\|^{1-\kappa_1}}+c_2\dfrac{\|x-x^*\|\cdot \|x^*-J_{\lambda A}(x-\lambda B(x))\|}{\|x-J_{\lambda A}(x-\lambda B(x))\|^{1-\kappa_2}}\right),
		\end{align*}
   for all $x\in H\setminus \{x^*\}$. 
		By the triangle inequality and Proposition \ref{tr1}, there exists \(\delta\in (0, 1)\) such that the following inequality
		\begin{equation} \label{ineq1}
			\|x-J_{\lambda A}(x-\lambda B(x))\|\leq \|x-x^*\|+\|J_{\lambda A}(x-\lambda B(x))-x^*\|\leq (1+\delta)\|x-x^*\|,
		\end{equation}
		holds for all \(x\in H\). Using the inverse triangle inequality we get   
  \begin{equation}\label{ineq2}
			\|x-J_{\lambda A}(x-\lambda B(x))\|\geq \|x-x^*\|-\|J_{\lambda A}(x-\lambda B(x))-x^*\|\geq (1-\delta)\|x-x^*\|,
		\end{equation}
		holds for all \(x\in H\).
  Combining inequalities \eqref{ineq1}, \eqref{ineq2} and Proposition \ref{tr1} we obtain
		\begin{align}\label{eq16} 
			\dot{V} \leq &-\left(\dfrac{c_1}{(1+\delta)^{1-\kappa_1}}\dfrac{\|x-x^*\|^2}{\|x-x^*\|^{1-\kappa_1}}+\dfrac{c_2}{(1+\delta)^{1-\kappa_2}}\dfrac{\|x-x^*\|^2}{\|x-x^*\|^{1-\kappa_2}}\right)\notag \\
			&+ \left(\dfrac{\delta c_1}{(1-\delta)^{1-\kappa_1}}\dfrac{\|x-x^*\|^2}{\|x-x^*\|^{1-\kappa_1}}+\dfrac{\delta c_2}{(1-\delta)^{1-\kappa_2}}\dfrac{\|x-x^*\|^2}{\|x-x^*\|^{1-\kappa_2}}\right) \notag \\
			=&-q(c_1, \kappa_1)\|x-x^*\|^{1+\kappa_1}-q(c_2, \kappa_2)\|x-x^*\|^{1+\kappa_2},
		\end{align}
		where \(q(c_i, \kappa_i):=\dfrac{c_i}{(1-\delta)^{1-\kappa_i}}\left(\left(\frac{1-\delta}{1+\delta}\right)^{1-\kappa_i}-\delta \right)\) for $i=1, 2.$
		By Lemma \ref{lm4}, there exists \(\eps(\delta)=\dfrac{\log(\delta)}{\log\frac{1-\delta}{1+\delta}}>0\) such that \(q(c_1, \kappa_1)>0\) for any \(\kappa_1\in (1-\eps(\delta), 1)\cap (0, 1)\) and \(q(c_2, \kappa_2)>0\) for any \(\kappa_2>1\). Hence, \eqref{eq16} implies that 
		\begin{equation}\label{eq17} 
			\dot{V} \leq -\left( p(c_1, \kappa_1)V^{\alpha(\kappa_1)}+p(c_2, \kappa_2)V^{\alpha(\kappa_2)} \right),
		\end{equation}
		here \(p(c_i, \kappa_i)=2^{\alpha(\kappa_i)}q(c_i, \kappa_i)\) and \(\alpha(\kappa_i)=\frac{1+\kappa_i}{2}\), for $i=1,2$. Note that \(p(c_1, \kappa_1)>0, \alpha(\kappa_1)\in (0.5, 1)\) for any \(\kappa_1\in (1-\eps(\delta), 1)\cap (0, 1)\) and \(p(c_2, \kappa_2)>0, \alpha(\kappa_2)>1\) for any \(\kappa_2>1\). Then the conclusion follows from Theorem \ref{lm1}.  
	\end{proof}
	
	\begin{remark}
		Assumption {\bf (A)} can be replaced by the following one: 
		\begin{itemize}
			\item [{\bf(A')}] The operator \(A\) is maximal strongly monotone with modulus \(\mu_A>0\) and \(B\) is \(\beta\)-cocoercive with \(\beta>0\). 
		\end{itemize}

		Then if assumption {\bf(A')} is valid, Proposition \ref{tr1}, Theorem \ref{tr2}  still hold whenever \(\lambda\in (0,2\beta)\).
		Indeed, with  assumption {\bf(A')}, the operator \(A+B\) is strongly monotone. Consequently, \(\mbox{zer}(A+B)\) is a singleton. Also, it can be proved that \(T:=J_{\lambda A}\circ (Id-\lambda B)\) is \(\delta\)-Lipschitz continuous with \(\delta\in (0,1)\) whenever \(\lambda\in (0,2\beta)\), (see, for example, the proof of Proposition 25.9 in \cite{BauschkeCombettes}). Therefore, Proposition \ref{tr1} still holds. The proof of Theorem \ref{tr2} under assumption {\bf(A')} is similar to the one above. 
		
	\end{remark}

	\section{Consistent discretization of the modified forward-backward dynamical system} \label{sec4}
	Continuous-time dynamical systems, exemplified by equation \eqref{newydynamicalsystem}, provide valuable insights for devising accelerated strategies to tackle inclusion problems. Nevertheless, in practice, discrete-time approaches are often applied with iterative methods. It's worth noting that while continuous-time dynamical systems exhibit fixed-time convergence behavior, this might not persist in discrete-time implementations. A consistent discretization scheme ensures that the convergence behavior of the continuous-time dynamical system is maintained in the discrete-time domain (refer to, for instance, \cite{POL}). In this section, we present a characterization of conditions conducive to achieving a consistent discretization of the fixed-time convergent modified forward-backward splitting dynamical system. This objective can be realized in a broader context of differential inclusions by leveraging the concepts outlined in \cite{BEN}.
	
	\begin{theorem}\cite{Garg21}\label{tr3} 
		Consider the following differential inclusion \begin{equation}\label{eq20} 
			\dot{x}\in \Xi(x)
		\end{equation}
		where \(\Xi: H \rightrightarrows H\) is an upper semi-continuous set-valued map satisfying that its values are non-empty, convex, compact, and  \(0\in \Xi(\overline{x})\) for some \(\overline{x}\in H\). Assume that there exists a positive definite, radially unbounded, locally Lipschitz continuous and regular function \(V: H\rightarrow H\) such that \(V(\overline{x})=0\) and 
		\begin{equation*}
			\sup \dot{V}(x)\leq -\left( r V(x)^{1-\frac{1}{\nu}}+s V(x)^{1+\frac{1}{\nu}}\right) 
		\end{equation*}
		for all \(x\in H\setminus\{\overline{x}\}\), with \(r, s>0\) and \(\nu>1\), here  
		\begin{equation*}
			\dot{V}(x)=\left\{ w\in \R: \exists u\in \Xi(x) \mbox{ such that } \langle z, u\rangle =w, \forall z\in \partial_cV(x)\right\}
		\end{equation*}
		and   \(\partial_c V(x)\) denotes Clarke's generalized gradient of the function \(V\) at the point \(x\in H\). Then, the equilibrium point \(\overline{x}\in H\) of \eqref{eq20} is fixed-time stable, with the settling-time function \(T\) satisfying 
		\begin{equation*}
			T(x(0)) \leq \dfrac{\nu\pi}{2\sqrt{rs}} 
		\end{equation*} 
		for any starting point  \(x(0)\in H\). 
	\end{theorem}
	We now examine  the forward-Euler discretization of \eqref{eq20}:
	\begin{equation}\label{eq24} 
		x_{n+1}\in x_n+\gamma \Xi(x_n)
	\end{equation}
	where \(\gamma >0\) is the time-step. 
	
	
	The theorem presented below characterizes the conditions under which a consistent discretization of a differential inclusion, featuring a fixed-time stable equilibrium point, is achieved.

	\begin{theorem}\cite{Garg21}\label{tr4} 
		Assume that the conditions of Theorem \ref{tr3} are fulfilled, and  the function \(V\) satisfies the following quadratic growth condition
		\begin{equation*}
			V(x)\geq \rho\|x-\overline{x}\|^2 
		\end{equation*}
		for every $x\in H$, where $\rho>0$ and 
  $\overline{x}$ 
  is the equilibrium point of \eqref{eq20}. 
  Then, for all \(x_0\in H\) and \(\eps>0\), there exists 
  {\color{blue} $\gamma^*>0$} such that for any \({\color{blue}\gamma}\in (0, {\color{blue}\gamma^*}]\), the following holds:
		\begin{equation*}
			\|x_n-\overline{x}\| <\begin{cases}
				\frac{1}{\sqrt{\rho}} \left( \sqrt{\frac{r}{s}}\tan\left(\frac{\pi}{2}-\frac{\sqrt{rs}}{\nu}{\color{blue}\gamma} n\right)\right)^{\frac{\nu}{2}}+\eps, &{\color{blue} \mbox{ if }\ n\leq n^*}\\
				\eps & \mbox{otherwise},
			\end{cases}
		\end{equation*}
		where \(x_n\) is a solution  of \eqref{eq24} starting from the point \(x_0\) {\color{blue} and $n^* = \left \lceil \dfrac{\nu\pi}{2\gamma\sqrt{rs}} \right \rceil$}.
	\end{theorem}
	
	Applying this result to the system \eqref{newydynamicalsystem} we get the following. 
	\begin{theorem}
		Consider the forward-Euler discretization of \eqref{newydynamicalsystem}:
		\begin{equation}\label{eq29} 
			x_{n+1}=x_n-\gamma \varphi(x_n)(x_n-J_{\lambda A}(x_n-\lambda B(x_n))),
		\end{equation}
		where  \(\varphi\) is given by \eqref{dnrho} with \(c_1, c_2>0, \kappa_1(\nu)=1-2/\nu\) and \(\kappa_2(\nu)=1+2/\nu\) where \(\nu\in (2, \infty)\), and \(\gamma>0\) is the time-step. Then for every \(x_0\in H\), every \(\eps>0\) and every \(\lambda>0\) satisfying conditions in assumption {\bf (A)}, there exist \(\nu>2, r, s>0\) and \(\gamma^*>0\) such that for any \(\gamma\in (0, \gamma^*]\), the following is fulfilled: 
		\begin{equation*}
			\|x_n-x^*\|< \begin{cases}
				\sqrt{2}\left( \sqrt{\frac{r}{s}}\tan\left(\frac{\pi}{2}-\frac{\sqrt{rs}}{\nu}\gamma n\right)\right)^{\frac{\nu}{2}}+\eps, & {\color{blue} \mbox{ if } n\leq n^*}\\ 
				\eps & \mbox{otherwise},
			\end{cases}
		\end{equation*}
		where {\color{blue}\(n^*=\left\lceil \dfrac{\nu\pi}{2\gamma \sqrt{rs}} \right\rceil\)} and \(x_n\) is a solution of \eqref{eq29} starting from the point \(x_0\) and \(x^*\in H\) is the unique solution of inclusion problem \eqref{inc1}. 
	\end{theorem}
	\begin{proof} We first note from the proof of Theorem \ref{tr2} that for any given \(\lambda\) satisfying the condition in {\bf (A)}, inequality \eqref{eq17} is valid for any \(\kappa_1(\nu)\in (1-\eps(\delta), 1)\cap (0, 1)\) with \(\eps(\delta)=\dfrac{\log(\delta)}{\log\left(\frac{1-\delta}{1+\delta}\right)}>0\) and \(\kappa_2(\nu)>1\). The former implies that \(\nu>\max\left\{2, \dfrac{2}{\eps(\delta)}\right\}\) and the latter is always fulfilled for any \(\nu>2\). Therefore, for any given \(\lambda\) as in assumption {\bf (A)}, inequality  \eqref{eq17} holds for all \(\nu>\max\left\{2, \dfrac{2}{\eps(\delta)}\right\}\) with \(p(c_1, \kappa_1(\nu))>0, p(c_2, \kappa_2(\nu))>0, \alpha\left( \kappa_1(\nu)\right) =1-\dfrac{1}{\nu}\) and \(\alpha\left(\kappa_2(\nu)\right) =1+\dfrac{1}{\nu}\). Thus, all conditions in Theorem \ref{tr4} are satisfied with the Lyapunov function $V(x)=\dfrac{1}{2}\|x-x^*\|^2$. The conclusion of the theorem is followed from Theorem \ref{tr4}. 
	\end{proof}
	
	
	This result shows that for every $\varepsilon>0$ the solution derived from applying forward-Euler discretization to \eqref{newydynamicalsystem} enters an $\varepsilon$-neighborhood of the solution to the associated inclusion problem within a fixed number of time steps, regardless of the initial conditions.

	\begin{remark}
		In the case \(B=0\), the inclusion problem \eqref{inc1} reduces to \(0\in A(x)\) and the the dynamical system \eqref{newydynamicalsystem} becomes 
		\begin{equation*}
			\dot{x}=-\varphi(x)(x-J_{\lambda A}(x)),
		\end{equation*} 
		where 
		\begin{equation*}
			\varphi(x):=\begin{cases}
				c_1\dfrac{1}{\|x-J_{\lambda A}(x)\|^{1-\kappa_1}}+c_2\dfrac{1}{\|x-J_{\lambda A}(x)\|^{1-\kappa_2}} & \mbox{ if } x\in H\setminus {\rm Fix}(J_{\lambda A})\\
				0 & \mbox{otherwise} 
			\end{cases}
		\end{equation*} 
		with \(c_1, c_2>0, \kappa_1\in (0, 1)\) and \(\kappa_2>1\).
		
		Then, all conclusions above still hold for solutions of the inclusion \(0\in A(x)\). 
	\end{remark}

	\section{Applications}\label{sec5}
	In this section, we apply the proposed fixed-times stable forward-backward splitting dynamical system \eqref{newydynamicalsystem} to address a range of optimization problems including constrained optimization problems (COPs), mixed variational inequalities (MVIs), and variational inequalities (VIs). 
	\subsection{Application to COP }
	Consider the constrained optimization problem (COP)
	\begin{equation}\label{cop1}
		\min_{x\in H} f(x)+g(x)
	\end{equation}
	where \(f: H\longrightarrow \R\) is continuous differential and convex, and \(g: H\longrightarrow \R\) is a proper, l.s.c. convex real value function. Observe that   \(g\) may not be differentiable and if \(g\equiv 0\), the constrained {\color{blue} problem} \eqref{cop1} become unconstrained optimization {\color{blue} one}. 
	
	Let \(T:=\nabla f+\partial g\). Then the problem of finding solutions of the COP \eqref{cop1} is equivalent to finding zero points of the operator \(T\), meaning that finding solutions of COP \eqref{cop1} is equivalent to finding solutions of inclusion problem \eqref{inc1} with $A=\partial g$ and $B=\nabla f$. 	
	Also, note that in this case 	
	$$J_{\lambda A}(x-\lambda B(x))=\mbox{prox}_{\lambda g}(x-\lambda \nabla f (x))$$
	and therefore, the dynamical system \eqref{newydynamicalsystem} becomes

	\begin{equation}\label{dDNMCOP}
		\dot{x}(t)= -\varphi(x)(x-\mbox{prox}_{\lambda g}(x-\lambda \nabla f(x))),
	\end{equation}
	where 
	\begin{equation*}
		\varphi(x):=\begin{cases}
			c_1\dfrac{1}{\|x-U(x)\|^{1-\kappa_1}} &+c_2\dfrac{1}{\|x-U(x)\|^{1-\kappa_2}} \mbox{ if } x\in H\setminus {\rm Fix}(U)\\
			0 & \mbox{otherwise} 
		\end{cases}
	\end{equation*} 
	with \(c_1, c_2>0, \kappa_1\in (0, 1)\) and \(\kappa_2>1\), $U:= \mbox{prox}_{\lambda g}\circ (Id-\lambda \nabla f)$.
	
	This is the proximal gradient dynamical system for COP studied in, for instance, \cite{ABBAS, HASS, JU} {\color{blue} where} \(\varphi(x)\equiv \varphi_0>0\) is a constant. 	As in \cite{JU5}, we make the following assumption on the function $f$.
	\begin{itemize}
		\item [{\bf (A1)}] The function $f$ is strongly monotone with a constant $\mu >0$ and $\nabla f$ is $L$-continuous with constant $L>0$, and $\lambda L^2<2\mu$.
	\end{itemize}

	Observe that when $f$ is strongly monotone with constant $\mu$,  $\nabla f$ is also strongly monotone with modulus $\mu$. Therefore, assumptions {\bf (A)} and {\bf (A')} hold when {\bf (A1)} is valid. 	By applying Theorem \ref{tr2} with \(A=\partial g\) and $B=\nabla f$ we obtain the fixed-time stable convergence for COP \eqref{cop1}.
	
	\begin{proposition}
		Let \(x^*\) be a solution to COP \eqref{cop1}. Assume that assumption {\bf (A1)} holds.  
		Then \(x^*\) is the fixed-times stable  equilibrium point of \eqref{dDNMCOP} with the settling time given as follows:
		\begin{equation*} 
			T(x(0)) \leq T_{\max} = \dfrac{1}{p_1(1-\alpha_1)}+\dfrac{1}{p_2(\alpha_2-1)}
		\end{equation*}
		where $p_1>0, p_2>0, \alpha_1\in (0.5,1),$ and $\alpha_2 >1$ are some constants.
		
		In addition, if take $\alpha_1=1-\frac{1}{2\xi}$ and \(\alpha_2=1+\frac{1}{2\xi}\) for some \(\xi>1\), then the following time estimate holds
		$$T(x(0))\leq T_{\max} =\dfrac{\pi \xi}{\sqrt{p_1p_2}}.$$
		\begin{proof}
			Its proof is followed by Theorem \ref{tr2} with noting that $J_{\lambda A}(x-\lambda B(x))=\mbox{prox}_{\lambda g}(x-\lambda \nabla f (x))$. 
		\end{proof}
	\end{proposition}
	\subsection{Application to MVIP}
	We now consider mixed variational inequality problem (MVIP): 
	\begin{equation}\label{mvip1}
		\mbox{Find \(x^*\in H\) such that \(\langle F(x^*), x-x^*\rangle +g(x^*)-g(x)\geq 0\) for all $x\in H$,}
	\end{equation}
	where $F: H\longrightarrow H$ is a vector-valued operator and $g: H\longrightarrow \R$ is a proper, l.s.c. convex function. 	These  MVIPs can be expressed in inclusion problem \eqref{inc1} with $B=F$ and $A = \partial g$. In this case, we have that 
	$$J_{\lambda A}(x-\lambda B(x))=\mbox{prox}_{\lambda g}(x-\lambda F(x)).$$
	Consequently, the dynamical system \eqref{newydynamicalsystem} becomes	
	\begin{equation}\label{dnm4}
		\dot{x}(t)=-\varphi(x)(x-\mbox{prox}_{\lambda g} (x-\lambda F(x))
	\end{equation}
	where 
	\begin{equation*}
		\varphi(x):=\begin{cases}
			c_1\dfrac{1}{\|x-\mbox{prox}_{\lambda g}(x-\lambda F(x))\|^{1-\kappa_1}} &+c_2\dfrac{1}{\|x-\mbox{prox}_{\lambda g}(x-\lambda F(x))\|^{1-\kappa_2}}\\
			& \mbox{ if } x\in H\setminus {\rm Fix}(\mbox{prox}_{\lambda g}\circ (Id-\lambda F))\\
			0 & \mbox{otherwise} 
		\end{cases}
	\end{equation*} 
	with \(c_1, c_2>0, \kappa_1\in (0, 1)\) and \(\kappa_2>1\).
	
	Then, it follows from Proposition \ref{pro4}  that $x^*\in H$ is a solution of MVIP \eqref{mvip1} if and only if it is an equilibrium point of the dynamical system \eqref{dnm4}.
	
	In the case where $g$ is convex,  $F$ is $\mu$-strongly monotone and $L$-Lipschitz  continuous  and $\lambda\in (0, 2\mu/L^2)$,   by Theorem \ref{tr2}, it holds that $x^*\in H$ is fixed-time stable solution of  MVIP \eqref{mvip1} with the following settling  time 
	\begin{equation*}  
		T(x(0)) \leq T_{\max} = \dfrac{1}{p_1(1-\alpha_1)}+\dfrac{1}{p_2(\alpha_2-1)},
	\end{equation*}
	where $p_1>0, p_2>0, \alpha_1\in (0.5,1),$ and $\alpha_2 >1$ are some constants.
	
	Moreover, if take $\alpha_1=1-\frac{1}{2\xi}$ and \(\alpha_2=1+\frac{1}{2\xi}\) for some \(\xi>1\), then the following settling time holds
	$$T(x(0))\leq T_{\max} =\dfrac{\pi \xi}{\sqrt{p_1p_2}}.$$
	\subsection{Application to VIP}
	We now turn back to the variational inequality problem:
	
	\begin{equation}\label{vip1}
		\mbox{Find \(x^*\in C\) such that \(\langle F(x^*), x-x^*\rangle  \geq 0, \quad \forall x\in C\),}
	\end{equation} 
 where \(C\) is a closed convex subset of \(H\), $F: C\longrightarrow H$ is an operator. 	We denote this problem by $\mbox{VIP}(F, C)$. Note that VIPs are a special case of MVIs when $g\equiv 0$. 
	
	The VIPs \eqref{vip1}  can be represented as the inclusion problem $0\in (A+B)(x)$ with $B=F$ and $A=N_C$. In this case  
 $$J_{\lambda A}(x-\lambda B(x)) =P_C(x-\lambda F(x)).$$	
	Consequently, the dynamical system \eqref{newydynamicalsystem} becomes
	\begin{equation}\label{dnm3}
		\dot{x}(t)=-\varphi(x)(x-P_C(x-\lambda F(x))),
	\end{equation}
	where 
	\begin{equation*}
		\varphi(x):=\begin{cases}
			c_1\dfrac{1}{\|x-U(x)\|^{1-\kappa_1}} &+c_2\dfrac{1}{\|x-U(x)\|^{1-\kappa_2}} \mbox{ if } x\in H\setminus {\rm Fix}(U)\\
			0 & \mbox{otherwise} 
		\end{cases}
	\end{equation*} 
	with \(c_1, c_2>0, \kappa_1\in (0, 1)\) and \(\kappa_2>1\), $U:= P_C\circ (Id-\lambda F)$.

	From the Proposition \ref{pro4} it holds that $x^*\in H$ is a solution of $\mbox{VIP}(F, C)$ if and only if \(x^*\) is an equilibrium point of the dynamical system \eqref{dnm3}.  In addition, if the operator \(F\) satisfies the following assumption: 
	\begin{itemize}
		\item [{\bf (A2)}] $F$ is strongly monotone with modulus \(\mu>0\) and Lipschitz continuous with constant \(L>0\)
	\end{itemize}	
	\noindent then, by applying Theorem \ref{tr2} to the operators \(A=N_C, B=F\)  we get the same result  as \cite [Theorem 2]{Garg21}.
 
	\begin{theorem}  	For every \(\lambda\in (0, \frac{2\mu}{L^2}) \), there exists \(\eps>0\) such that the solution \(x^*\in H\) of the variational inequality  problem $\mbox{VIP}(F, C) $ is a fixed-time stable equilibrium point of \eqref{dnm3} for any \(\kappa_1\in (1-\eps, 1)\cap(0, 1)\) and \(\kappa_2>1\).  Also, the following settling time holds
		\begin{equation*}  
			T(x(0)) \leq T_{\max} = \dfrac{1}{p_1(1-\alpha_1)}+\dfrac{1}{p_2(\alpha_2-1)}
		\end{equation*}
		where $p_1>0, p_2>0, \alpha_1\in (0.5,1),$ and $\alpha_2 >1$ are some constants.
		
		Moreover, if $\alpha_1=1-\frac{1}{2\xi}$ and \(\alpha_2=1+\frac{1}{2\xi}\) for some \(\xi>1\). Then the following settling time holds
		$$T(x(0))\leq T_{\max} =\dfrac{\pi \xi}{\sqrt{p_1p_2}}.$$
	\end{theorem}

\section{Numerical Illustration}\label{num}
We present in this section some numerical examples to illustrate the theoretical results obtained on the previous sections. The time-continuous dynamical systems
are solved by Runge–Kutta adaptive method (ode45 in MATLAB) on the time interval $[0,T_{max}]$.\\

\begin{figure}
	\begin{centering}
		\includegraphics[width=\textwidth]{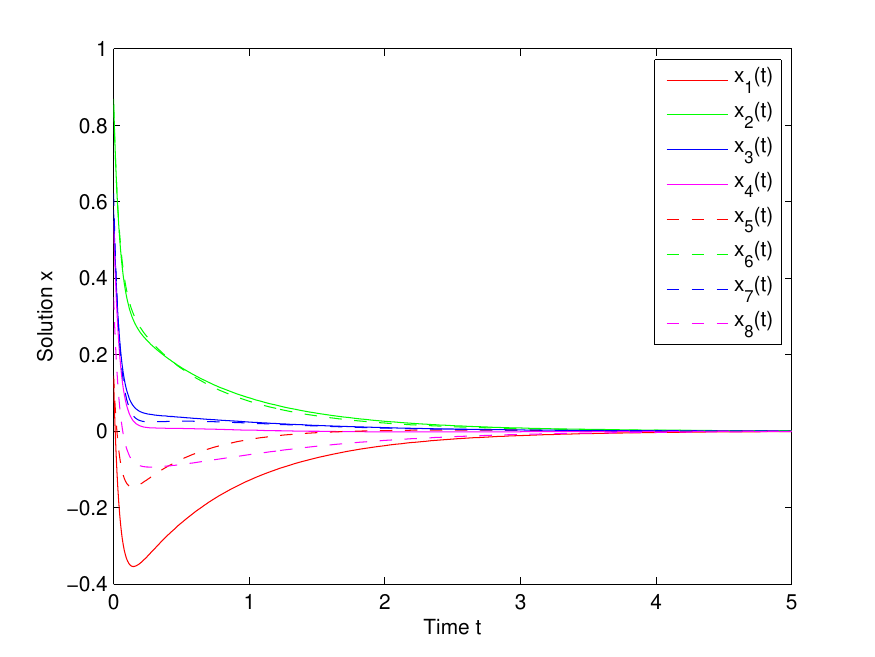}
		\par\end{centering}
	\protect\caption{Transient responses of dynamic system \eqref{newydynamicalsystem} for Example 1.
		\label{fig1}}
\end{figure}

\begin{figure}
	\begin{centering}
		\includegraphics[width=\textwidth]{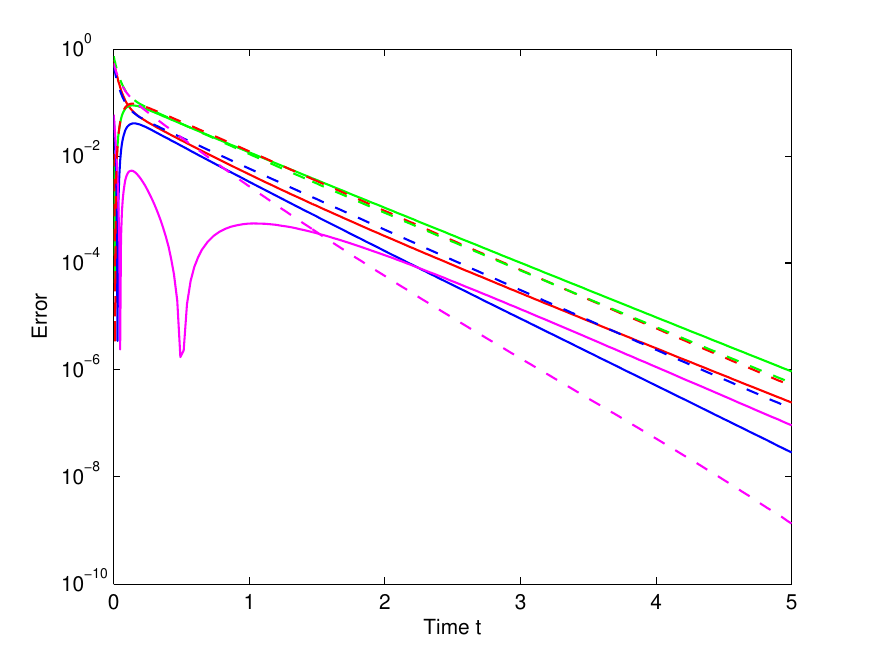}
		\par\end{centering}
	\protect\caption{Error convergence
responses of dynamic system \eqref{newydynamicalsystem} for Example 1.
		\label{fig2}}
\end{figure}

{\bf Example 1:}
We consider a monotone inclusion of the form \eqref{inc1} with
$$
A = I, \quad B = M^TM,
$$
where $I$ is the identity matrix of size $n$ and $M$ is an arbitrary matrix of order $m \times n$.
It is clear that the monotone inclusion $
0 \in A(x) +B(x)$ has a unique solution $x^*=0$. 
In our experiments, we choose $m=10, n=8, c_1=20, c_2=200, \kappa_1=0.99, \kappa_2=1.01$ and $\lambda = 0.005$. By applying the dynamic system \eqref{newydynamicalsystem} with $T_{max}=5$, Figure \ref{fig1} represents the solution of the above monotone inclusion which converges to $x^*$.  Figure \ref{fig2} shows the error convergence $\|x(t)-x^*\|^2$ of the solution $x(t).$\\

{\bf Example 2:}  We consider a COP which are widely used in statistics, machine learning, and signal
processing:
$$ \min_{x\in R^n} f(x)+ g(x)$$ where $f(x)=\dfrac{1}{2}\|Qx-a\|^2, g(x) = \alpha \|x\|_1, Q$ is {\color{blue} an  $m\times n$ matrix,} $ a \in R^m$ is the response vector and $\alpha > 0$ is the regularization parameter.
Let $U=\nabla f + \partial g$, then the solution of the COP is zero points of the operator $U$. In other words, solving the COP	is equivalent to solving the inclusion problem (\ref{inc1}) with $A(x)=\partial g(x)$ and $B(x)=\nabla f(x) = Q^T(Qx-a)$. \\

\begin{figure}
	\begin{centering}
		\includegraphics[width=\textwidth]{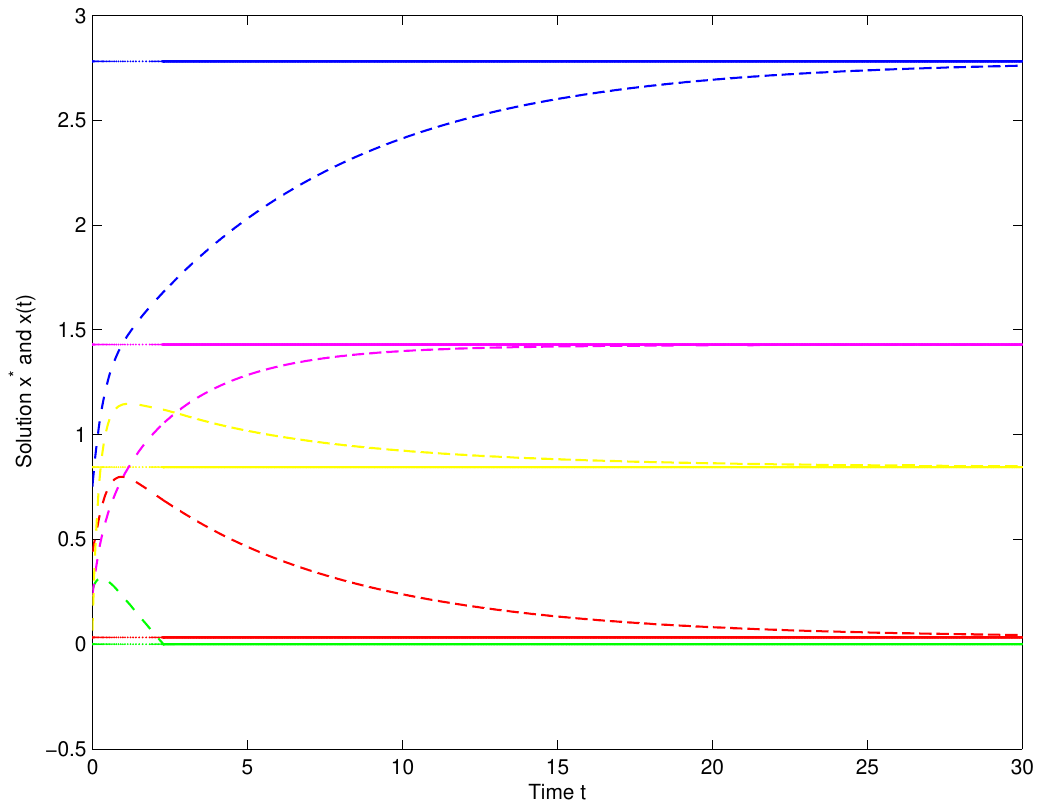}
		\par\end{centering}
	\protect\caption{Transient responses of dynamic system \eqref{dDNMCOP} and solution $x^*$ of Example 2.
		\label{fig3}}
\end{figure}

\begin{figure}
	\begin{centering}
		\includegraphics[width=\textwidth]{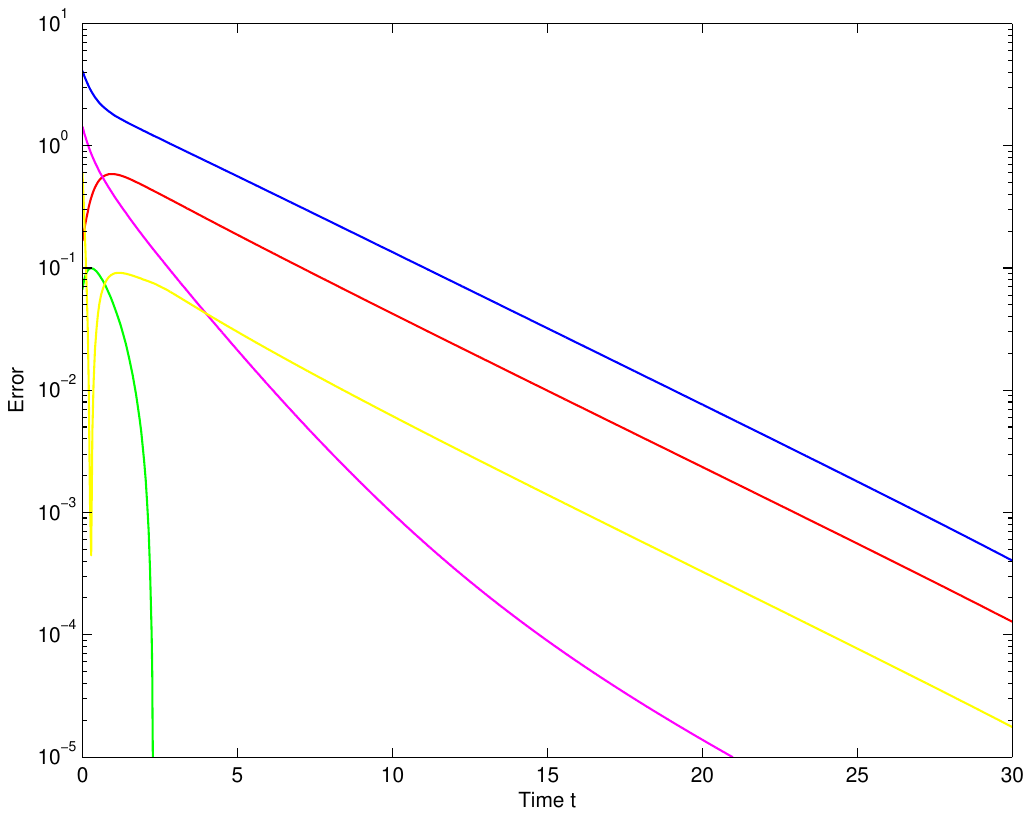}
		\par\end{centering}
	\protect\caption{Error convergence
responses of dynamic system \eqref{dDNMCOP} for Example 2.
		\label{fig4}}
\end{figure}

In the illustration, we choose  
$$
Q = \left[ \begin{array}{ccccc}5 & 0 & 2 & 0 & -1 \\-1 & 1 & 1 & 0 & 6 \\ 2 & 0 & 2 & 1 & -1 \\ 0 & 0 & 0 & 2 & 3 \end{array}\right],\quad
a = \left[\begin{array}{r}5\\8\\7\\6 \end{array}\right]; \quad \alpha = 2. 
$$
It is clear 
that $\nabla f$ is L-continuous and strongly monotone. Dynamical system (\ref{dDNMCOP}) with $c_1=20, c_2=200, \kappa_1=0.98, \kappa_2=1.02, \lambda = 0.005$ and $T_{max} =30$ creates the solution $x(t)$ which converges to the solution $x^*$ of the COP as shown in  Figure \ref{fig3}. Notice that with {\color{blue} a} randomly initial point, the {\bf solid lines} are for the solution $x^*$ obtained by using "fmincon" function in MATLAB and {\bf dashed lines} are for the solution $x(t)$. The error convergence $\|x(t)-x^*\|^2$ is also illustrated in  Figure \ref{fig4}. \\

\begin{figure}
	\begin{centering}
		\includegraphics[width=\textwidth]{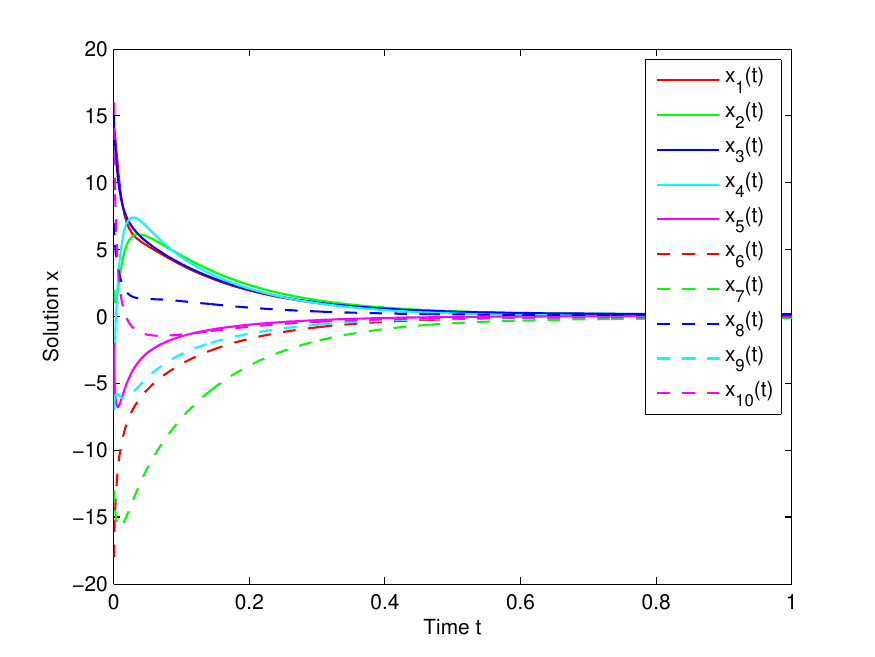}
		\par\end{centering}
	\protect\caption{Transient responses of dynamic system \eqref{dnm3} for Example 3.}
		\label{fig5}
\end{figure}

{\bf Example 3:} In this experiment, we consider a VIP: find $x^* \in C$ such that $$\langle F(x^*),x-x^* \rangle \geq 0, \quad \forall x \in C,$$ 
where $$C=\{x \in R^n: -20 \leq x_i \leq 20,\quad i=1,2,...n\}.$$ 
We choose 
$F(x) = Gx+u$ and $u \in R^n$ is a given vector with
$$ G = M^TM + S + T, $$
where $M$ is a square matrix of size $n$, $S$ is a skew-symmetric matrix of size $n$ and $T$ is a diagonal matrix of size $n$ with its diagonal entries being positive. 
The VIP can be represented as the inclusion problem $0 \in (A+B)(x)$ with $B = F$ and $A = N_C$ that assumption {\bf (A)} is satisfied. For the experiment, all entries of $u, M,$ and $S$ are generated randomly in $(-5,5)$, while entries of $T$ are generated randomly in (0,5), and $n =10$. 
The convergence of the solution $x(t)$ generated by the dynamic system \eqref{dnm3} is shown in the Figure \ref{fig5} with parameter set: $c_1=100, c_2=100, \kappa_1 = 0.9, \kappa_2 = 1.1, \lambda = 0.01$ and $T_{max} = 1 $.

\begin{figure}
	\begin{centering}
		\includegraphics[width=\textwidth]
  {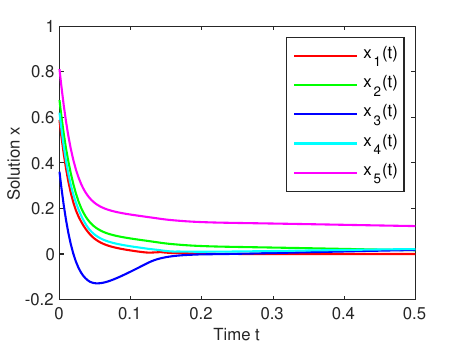}
		\par\end{centering}
	\protect\caption{Transient responses of dynamic system (\ref{dDNMCOP}) for Example 4.}
		\label{fig6}
\end{figure}

{\bf Example 4:} We consider the following elastic-net logistic regression problem with an $L^1$-regularization term:
$$
\min_{x\in \R^5}\sum_{i=1}^{100}\{\log(1+\exp(-p_iq_i^Tx))+\beta\|x\|_1\}
$$
where $p_i \in \{-1,1\}, q_i \in \R^5, i = 1,2,3,...100$ are chosen randomly, $\beta > 0.$ The above elastic-net logistic regression problem can be solved by the dynamic system (\ref{dDNMCOP}), where $f(x) = \sum_{i=1}^{100}{\log(1+\exp(-p_iq_i^Tx))}$ and $g(x)=\beta\|x\|_1$. The parameter set chosen for numerical experimentation is: $\beta = 2.5, c_1=10, c_2=10, \kappa_1 = 0.9, \kappa_2 = 1.1, \lambda = 0.01$.  Figure \ref{fig6} demonstrates the convergence of the solution $x(t)$ 
 obtained using the dynamic system (\ref{dDNMCOP}) 
 with $T_{max} =0.5$
 to the solution $x^*$.
\\

\newpage

	\section{Conclusion}\label{sec6}
	
	This paper provides a modified forward-backward splitting dynamical system, ensuring the existence, uniqueness, and convergence of its solution to the unique solution of the inclusion problems within a fixed-time context. This achievement is attained under a generalized monotonicity of involved operators. Additionally, the paper demonstrates that the forward-Euler discretization of the modified dynamical system serves as a consistent discretization method. Finally, the paper presents some applications by treating MVIs, VIs, and COPs as special cases of inclusion problems.  Future research directions include exploring the fixed-time stability of the modified forward-backward splitting dynamical system in the broader context of infinite-dimensional Hilbert or Banach spaces and investigating its stability under relaxed assumptions akin to those in prior studies.

\section*{Acknowledgement} 	
The authors thank the anonymous referee for their helpful comments.

\section*{Declarations} 
{\bf Conflict of interest} The authors declare no competing interests.

\end{document}